\documentclass[amssymb, 11pt]{amsart}
\usepackage{latexsym}

\newlength{\standardunitlength}
\newtheorem{prop}{Proposition}[section]

\newtheorem{lemma}[prop]{Lemma}
\newtheorem{cor}[prop]{Corollary}
\newtheorem{theorem}[prop]{Theorem}

\begin{document}

\begin{center}
{\bf Stein's Method, Jack Measure, and the Metropolis Algorithm}
\end{center}

\begin{center}
{\bf Running head: Stein's Method and Jack Measure}
\end{center}

\begin{center}
Written 11/13/03; Referee suggestions implemented on 6/30/04.
\end{center}

\begin{center}
By Jason Fulman
\end{center}

\begin{center}
University of Pittsburgh
\end{center}

\begin{center}
Department of Mathematics
\end{center}

\begin{center}
301 Thackeray Hall
\end{center}

\begin{center}
Pittsburgh, PA 15260
\end{center}

\begin{center}
Email: fulman@math.pitt.edu
\end{center}

{\bf Abstract}: The one parameter family of Jack$_{\alpha}$ measures
on partitions is an important discrete analog of Dyson's $\beta$
ensembles of random matrix theory. Except for special values of
$\alpha=1/2,1,2$ which have group theoretic interpretations, the
Jack$_{\alpha}$ measure has been difficult if not intractable to
analyze. This paper proves a central limit theorem (with an error
term) for Jack$_{\alpha}$ measure which works for arbitrary values of
$\alpha$. For $\alpha=1$ we recover a known central limit theorem on
the distribution of character ratios of random representations of the
symmetric group on transpositions. The case $\alpha=2$ gives a new
central limit theorem for random spherical functions of a Gelfand
pair (or equivalently for the spectrum of a natural random walk on perfect matchings in the complete graph). The proof uses Stein's method and has interesting combinatorial
ingredients: an intruiging construction of an exchangeable pair,
properties of Jack polynomials, and work of Hanlon relating Jack
polynomials to the Metropolis algorithm.

\begin{center}
\end{center}

\begin{center}
 \end{center}

\newpage

\section{Introduction}

	The purpose of this paper is to give a new approach to
studying a certain probability measure on the set of all partitions of
size $n$, known as Jack$_{\alpha}$ measure. Here $\alpha>0$, and this
measure chooses a partition $\lambda$ of size $n$ with probability \[
\frac{\alpha^n n!}{\prod_{s \in \lambda} (\alpha a(s) + l(s) +1)
(\alpha a(s) + l(s) + \alpha)}, \] where the product is over all boxes
in the partition. Here $a(s)$ denotes the number of boxes in the same
row of $s$ and to the right of $s$ (the ``arm'' of s) and $l(s)$
denotes the number of boxes in the same column of $s$ and below $s$
(the ``leg'' of s). For example the partition of 5 below \[
\begin{array}{c c c} \framebox{}& \framebox{}& \framebox{} \\
\framebox{}& \framebox{}& \end{array}, \] would have Jack$_{\alpha}$
measure \[ \frac{60 \alpha^2}{(2 \alpha+2)(3 \alpha+1) (\alpha+2)(2
\alpha+1)(\alpha+1)}.\]

Before proceeding, it should be mentioned that there is significant
interest in the study of statistical properties of Jack$_{\alpha}$
measure when $\alpha$ is fixed. The case $\alpha=1$ corresponds to the
Plancherel measure of the symmetric group, which is now well
understood due to numerous results in the past few years. The surveys
\cite{AD},\cite{De}, \cite{O2} and the seminal papers
\cite{BOO},\cite{J},\cite{O1} indicate how the Plancherel measure of
the symmetric group is a discrete analog of random matrix theory, and
describe its importance in representation theory and geometry. The
case $\alpha=2$ corresponds to the Gelfand pair $(S_{2n},H_{2n})$
where $S_{2n}$ is a symmetric group and $H_{2n}$ is the
hyperoctahedral group of size $2^nn!$. When $\alpha=\frac{1}{2}$, Jack
polynomials arise in the study of the Gelfand pair $(GL(n,H),U(n,H))$
where $H$ denotes the division ring of quaternions and $GL,U$ denote
general linear and unitary group. The paper \cite{O2} emphasizes that
the study of Jack$_{\alpha}$ measure is an important open problem,
about which relatively little is known \cite{BO1}. It is a discrete
analog of Dyson's $\beta$ ensembles, which are tractable for the three
values $\beta=1,2,4$. In particular, the correlation functions of
Jack$_{\alpha}$ measure are not known, so the traditional techniques
for studying discrete analogs of random matrix theory are not
obviously applicable.

	In the current paper we study Jack$_{\alpha}$ measure using a
remarkable probability technique known as Stein's method. Although
Stein's method can be quite hard to work with, there are some problems
where it seems to be the only option available (see \cite{RR} for such
an example involving the antivoter model). Good surveys of Stein's
method (two of them books) are \cite{AGG}, \cite{BHJ}, \cite{St1},
\cite{St2}.

	The current paper is a continuation of \cite{F1}, which applied Stein's method to the study of Plancherel measure of the symmetric group $S_n$. Let $\chi^{\lambda}_{(2,1^{n-2})}$ denote the character of the irreducible representation of $S_n$ parameterized by $\lambda$ on the conjugacy class of transpositions. Let $dim(\lambda)$ denote the dimension of the irreducible representation parameterized by $\lambda$. Letting $P_{\alpha}$ denote the probability of an event under Jack$_{\alpha}$ measure (so that $P_{1}$ corresponds to Plancherel measure), the following central limit theorem was proved:
\begin{theorem} \label{mainold} (\cite{F1}) For $n \geq 2$ and all real $x_0$, \[ |P_{1} \left( \frac{n-1}{\sqrt{2}} \frac{\chi^{\lambda}_{(2,1^{n-2})}}{dim(\lambda)} \leq
x_0 \right) - \frac{1}{\sqrt{2 \pi}} \int_{-\infty}^{x_0}
e^{-\frac{x^2}{2}} dx| \leq 40.  1 n^{-1/4}. \]
\end{theorem} This result sharpened earlier work of Kerov \cite{K1} (see \cite{IO} for a detailed exposition of Kerov's argument) and Hora \cite{Ho}, who
 both obtained a central limit theorem by the method of moments, but
 with no error bound. We remark that statistical properties of the
 quantity $\frac{\chi^{\lambda}_{(2,1^{n-2})}}{dim(\lambda)}$ (also called a
 character ratio) have important applications to random walk \cite{DS}
 and to the moduli space of curves \cite{EO}.

	The main result of the current paper is the following
deformation of Theorem \ref{mainold}. To state it one needs some
notation about partitions. Let $\lambda$ be a partition of some
non-negative integer $|\lambda|$ into integer parts $\lambda_1 \geq
\lambda_2 \geq \cdots \geq 0$. The symbol $m_i(\lambda)$ will denote
the number of parts of $\lambda$ of size $i$. Let $l(\lambda)$ denote
$\sum_{i \geq 1} m_i(\lambda)$, the number of parts of $\lambda$. Let
$n(\lambda)$ be the quantity $\sum_{i \geq 1} (i-1) \lambda_i$. One
defines $\lambda'$ to be the partition dual to $\lambda$ in the sense
that $\lambda_i' = m_i(\lambda) + m_{i+1}(\lambda) +
\cdots$. Geometrically this corresponds to flipping the diagram of
$\lambda$.

\begin{theorem} \label{maintheorem} Suppose that $\alpha \geq 1$. Let $W_{\alpha}(\lambda) = \frac{\alpha n(\lambda')-n(\lambda)}{\sqrt{\alpha {n \choose 2}}}$. For $n \geq 2$ and all real $x_0$, \[ |P_{\alpha} \left( W_{\alpha} \leq
x_0 \right) - \frac{1}{\sqrt{2 \pi}} \int_{-\infty}^{x_0}
e^{-\frac{x^2}{2}} dx| \leq A_{\alpha} n^{-1/4} \] where $A_{\alpha}$ depends on $\alpha$ but not on $n$.
\end{theorem}

	Note that the assumption that $\alpha \geq 1$ is merely for
convenience. Indeed, from the definition of Jack measure it is clear
that the Jack$_{\alpha}$ probability of $\lambda$ is equal to the
Jack$_{\frac{1}{\alpha}}$ probability of $\lambda'$. From this one
concludes that the Jack$_{\alpha}$ probability that $W_{\alpha}=w$ is
equal to the Jack$_{\frac{1}{\alpha}}$ probability that
$W_{\frac{1}{\alpha}}=-w$, so that a central limit theorem holds for
$\alpha$ if and only if it holds for $\frac{1}{\alpha}$.

	We conjecture that the convergence
rate upper bound in Theorem \ref{maintheorem} can be improved to a universal constant multiplied by the maximum of $\frac{1}{\sqrt{n}}$ and $\frac{\sqrt{\alpha}}{n}$. In fact the third moment of
$W_{\alpha}$ is $\frac{\alpha-1}{\sqrt{\alpha {n \choose 2}}}$ (see
Corollary \ref{exvar}), so certainly $\frac{\sqrt{\alpha}}{n}
\rightarrow 0$ is necessary for $W_{\alpha}$ to be asymptotically
normal. Of course typically one is interested in $\alpha$ fixed, as
$\alpha$ is a parameter which represents the symmetries of the system.

A result of Frobenius \cite{Fr} is that \[
 \frac{\chi^{\lambda}_{(2,1^{n-2})}}{dim(\lambda)} =
 \frac{n(\lambda')-n(\lambda)}{{n \choose 2}}. \] Hence Theorem
 \ref{maintheorem} is a generalization of Theorem \ref{mainold} in the
 case $\alpha=1$. It is also of group theoretic interest in the case
 $\alpha=2$. By
 page 410 of \cite{Mac} one sees for the $\alpha=2$ case that
 $\frac{2n(\lambda')-n(\lambda)}{2{n \choose 2}}$ is the value of a
 spherical function corresponding to the Gelfand pair
 $(S_{2n},H_{2n})$, where $H_{2n}$ is the hyperoctahedral group of
 size $2^n n!$. Moreover when $\alpha=2$, Theorem \ref{maintheorem} gives a
central limit theorem for the spectrum of a natural random walk on
perfect matchings of the complete graph. For a definition and analysis
of the convergence rate of this random walk on matchings, see
\cite{DHol}, where it was studied in connection with phylogenetic
trees. Note that their Corollary 1 shows that the eigenvalues of that
random walk are indexed by partitions $\lambda$ of $n$, and are
$\frac{W_2(\lambda)}{\sqrt{n(n-1)}}$, occurring with multiplicity
proportional to the Jack$_2$ measure on $\lambda$.

	Next we make some remarks about the proof of Theorem
\ref{maintheorem}. The argument is not a straightforward modification
of arguments used in \cite{F1}, and requires new
ideas. The reason for this is that for general $\alpha$ the
Jack$_{\alpha}$ measure does not have a known interpretation in terms
of representation theory of finite groups. Hence the proof of
\cite{F1}, which used concepts such as induction and restriction of
characters, can not be applied. There is another fundamental
difference between the case of Plancherel measure and Jack$_{\alpha}$
measure. In the Plancherel case the argument of \cite{F1} can be
pushed through to conjugacy classes other than transpositions, but the
same is not clearly so for the Jack$_{\alpha}$ case. This is because
the Jack$_{\alpha}$ case uses connections between Jack polynomials and
the Metropolis algorithm (due to Hanlon \cite{Ha} and to be reviewed
in Section \ref{metropjack}) and it is not clear that these
connections work for classes other than transpositions.

	Theorem \ref{maintheorem} will be a consequence of the
	following bound of Stein. Recall that if $W,W^*$ are random
	variables, they are called exchangeable if for all $w_1,w_2$,
	$P(W=w_1,W^*=w_2)$ is equal to $P(W=w_2,W^*=w_1)$. The
	notation $E^W(\cdot)$ means the expected value given $W$. Note from 
\cite{St1} that there are minor variations on Theorem \ref{steinbound} (and thus 
for Theorem \ref{maintheorem}) for $h(W)$ where $h$ is a bounded continuous 
function with bounded piecewise continuous derivative. For simplicity we only 
state the result when $h$ is the indicator function of an interval. 

\begin{theorem} \label{steinbound} (\cite{St1}) Let $(W,W^*)$ be an exchangeable pair
of real random variables such that $E^W(W^*) = (1-\tau)W$ with
$0<\tau<1$. Then for all real $x_0$, \begin{eqnarray*} & & |P(W \leq x_0) -
\frac{1}{\sqrt{2 \pi}} \int_{-\infty}^{x_0} e^{-\frac{x^2}{2}} dx|\\
& \leq & 2 \sqrt{E[1-\frac{1}{2 \tau}E^W (W^*-W)^2]^2} + (2 \pi)^{-
\frac{1}{4}} \sqrt{\frac{1}{\tau} E|W^*-W|^3}. \end{eqnarray*} \end{theorem}

	In order to apply Theorem \ref{steinbound} to study a
	statistic $W$, one clearly needs an exchangeable pair
	$(W,W^*)$ such that $E^W(W^*) = (1-\tau)W$. A Markov chain $K$
	(with chance of going from $x$ to $y$ denoted by $K(x,y)$) on
	a finite set $X$ is called reversible with respect to a
	probability distribution $\pi$ if $\pi(x) K(x,y) = \pi(y)
	K(y,x)$. This condition implies that $\pi$ is a stationary
	distribution for $K$. The idea is to use a reversible Markov
	chain on the set of partitions of size $n$ whose stationary
	distribution is Jack$_{\alpha}$ measure, to let $\lambda^*$ be
	obtained from $\lambda$ by one step in the chain where
	$\lambda$ is sampled from $\pi$, and then set
	$(W,W^*)=(W(\lambda),W(\lambda^*))$. A main contribution of
	this paper is the construction and analysis of an exchangeable
	pair which is useful for Stein's method.

Section \ref{revisit} revisits and generalizes the construction of an
exchangeable pair for Plancherel measure of the symmetric group. We
give a connection between harmonic functions on Bratelli diagrams and
decomposition of tensor products and extend some results in
\cite{F2}. Section \ref{jackprop} reviews necessary facts about Jack
polynomials. Motivated by the discussion in Section \ref{revisit},
Section \ref{exchange} constructs an exchangeable pair
$(W_{\alpha},W^*_{\alpha})$ to be used in the proof of Theorem
\ref{maintheorem}. The combinatorics in this section is quite
interesting. Section \ref{metropjack} recalls needed work of Hanlon
\cite{Ha} relating Jack polynomials to the Metropolis
algorithm. Section \ref{CLT} combines the ingredients of the previous
sections to prove Theorem \ref{maintheorem}.

	To close the introduction, we mention some follow up work to
this paper. The paper \cite{F3} sharpens the bound in Theorem
\ref{maintheorem} using martingale theory. The forthcoming paper
\cite{CF} extends the approach of this paper to other Gelfand pairs
(where the limit need not be a Gaussian law).

\section{Plancherel Measure Revisited} \label{revisit} 

	To begin, we revisit the construction of an exchangeable pair $(W,W')$ for the special case $\alpha=1$, corresponding to Plancherel measure, which was studied in \cite{F1}. In doing so we clarify and generalize some of the results there and in \cite{F2}. This will be very helpful for treating the case of general $\alpha$.

	As mentioned in the introduction, to construct an exchangeable
pair $(W,W^*)$ with respect to a probability measure $\pi$ on a finite
set $X$, it is enough to construct a Markov chain on $X$ which is
reversible with respect to $\pi$. Indeed, choosing $x$ from $\pi$ and letting $x^*$ be obtained
from $x$ by one step of the chain, it follows that $(W,W^*) :=
(W(x),W(x^*))$ is an exchangeable pair. Of course one wants to
construct the Markov chain in such a way that the exchangeable pair is
useful for Stein's method, and more precisely useful for Theorem
\ref{steinbound}.

\subsection{Known Constructions}

To start we consider the situation for an arbitrary finite group $G$. Let $Irr(G)$ denote the set of irreducible representations of $G$. Then the Plancherel measure on $Irr(G)$ chooses a representation $\lambda$ with probability $\frac{dim(\lambda)^2}{|G|}$, where $dim(\lambda)$ denotes the dimension of $\lambda$. In \cite{F2} we constructed a Markov chain $M_H$ on $Irr(G)$ which is reversible with respect to Plancherel measure. To define this Markov chain, one first fixes a subgroup $H$ of $G$. For $\tau \in Irr(H)$ and $\rho \in Irr(G)$, we let $\kappa(\tau,\rho)$ denote the multiplicity of $\rho$ in the character of $G$ obtained by inducing $\tau$ from $H$ (by Frobenius reciprocity, this is also equal to the multiplicity of $\tau$ in the character of $H$ obtaining by restricting $\rho$). Then \cite{F2} defined the transition probability $M_H(\lambda,\rho)$ of moving from a representation $\lambda$ to a representation $\rho$ by \[ \frac{|H|}{|G|} \frac{dim(\rho)}{dim(\lambda)} \sum_{\tau \in Irr(H)} \kappa(\tau,\lambda) \kappa(\tau,\rho). \] It was proved there that these transition probabilities sum to one, and that the Markov chain with transition mechanism $M_H$ is indeed reversible with respect to the Plancherel measure of $G$.

	For arbitrary groups, this construction can be recast in terms of harmonic functions on Bratelli diagrams. We recommend \cite{K2} or \cite{BO2} for an introduction to this subject. One starts with a
Bratteli diagram; that is an oriented graded graph $\Gamma= \cup_{n
\geq 0} \Gamma_n$ such that

\begin{enumerate}
\item $\Gamma_0$ is a single vertex $\emptyset$.
\item If the starting vertex of an edge is in $\Gamma_i$, then its end
vertex is in $\Gamma_{i+1}$.
\item Every vertex has at least one outgoing edge.
\item All $\Gamma_i$ are finite.
\end{enumerate}

For two vertices $\lambda, \Lambda \in \Gamma$, one writes $\lambda \nearrow \Lambda$ if there is an edge from $\lambda$ to
$\Lambda$. Part of the underlying data is a multiplicity function
$\kappa(\lambda,\Lambda)$. Letting the weight of a path in $\Gamma$ be
the product of the multiplicities of its edges, one defines the
dimension $dim(\Lambda)$ of a vertex $\Lambda$ to be the sum of the
weights over all minimal length paths from $\emptyset$ to $\Lambda$. Given a Bratteli
diagram with a multiplicity function, one calls a function $\phi$ {\it
harmonic} if $\phi(0)=1$, $\phi(\lambda) \geq 0$ for all $\lambda \in
\Gamma$, and \[ \phi(\lambda) = \sum_{\Lambda: \lambda \nearrow
\Lambda} \kappa(\lambda,\Lambda) \phi(\Lambda).\] An equivalent concept is that of coherent probability distributions. Namely a set
$\{M_n\}$ of probability distributions $M_n$ on $\Gamma_n$ is called
{\it coherent} if \[ M_{n-1}(\lambda) = \sum_{\Lambda: \lambda
\nearrow \Lambda} \frac{dim(\lambda)
\kappa(\lambda,\Lambda)}{dim(\Lambda)} M_{n}(\Lambda).\] The formula showing the concepts to be equivalent is $\phi(\lambda) =
\frac{M_n(\lambda)}{dim(\lambda)}$. Note that in this setting there is a natural transition mechanism for moving up or down a step in the Bratelli diagram. Namely the chance of moving from $\lambda$ to $\Lambda$ is $\frac{\kappa(\lambda,\Lambda) M_n(\Lambda) dim(\lambda)}{M_{n-1}(\lambda) dim(\Lambda)}$, and the chance of moving from $\Lambda$ to $\lambda$ is $\frac{dim(\lambda) \kappa(\lambda,\Lambda)}{dim(\Lambda)}$.

	Let $H_0=\{id\} \subseteq H_1 \subseteq \cdots \subseteq
	H_n=G$ be a tower of subgroups of $G$. Consider the Bratelli
	diagram whose $j$th level consists of irreducible
	representations of $H_j$, with edge multiplicity given by
	$\kappa(\tau,\lambda)$ as in the first paragraph of this
	subsection. It is proved in \cite{F2} that the Plancherel
	measures of the groups form a coherent family of probability
	distributions (this was known for the symmetric group
	\cite{K1}). Moreover it was shown that if one transitions from
	level n to level n-1, and then from level n-1 to level n, that
	the resulting Markov chain on irreducible representations of
	$H_n$ is exactly the chain $M_{H_{n-1}}$.

\subsection{New Construction}

Next we give a new Markov chain $L_{\eta}$ on the set of irreducible representations of $G$ which is reversible with respect to Plancherel measure, and which generalizes the chain $M_H$. First fix $\eta$, any representation (not necessarily irreducible) of $G$ whose character is real valued. Let $<\phi,\psi>$ be the usual inner product on class functions of $G$ defined as $\frac{1}{|G|} \sum_{g \in G} \phi(g) \overline{\psi(g)}$. Then the probability that the chain $L_{\eta}$ transitions from $\lambda$ to $\rho$ is \[ \frac{dim(\rho)}{dim(\eta) dim(\lambda)} <\chi^{\rho},\chi^{\eta} \chi^{\lambda}>.\] Note that this is nonnegative because $<\chi^{\rho},\chi^{\eta} \chi^{\lambda}>$ is the multiplicity of $\rho$ in the tensor product of $\eta$ and $\lambda$. 

\begin{lemma} \label{usereal} Let $\eta$ be a representation of a finite group $G$ whose character is real valued. Then the transition probabilities of $L_{\eta}$ sum to 1, and the Markov chain $L_{\eta}$ is reversible with respect to the Plancherel measure of $G$. \end{lemma}

\begin{proof} To see that the transition probabilities do as claimed sum to 1, observe
 that $\sum_{\rho} dim(\rho) \chi^{\rho}$ is the character of the
 regular representation of $G$, so takes value $|G|$ at the identity element and 0 elsewhere. The reversibility assertion uses the
 fact that $<\chi^{\rho},\chi^{\eta} \chi^{\lambda}>$ is equal to
 $<\chi^{\eta} \chi^{\rho},\chi^{\lambda}>$, which is true since
 $\chi^{\eta}$ is real valued. \end{proof}

	We remark that the second part of Lemma \ref{usereal} needs $\chi^{\eta}$ to be real
	valued. An instructive counterexample when $\chi^{\eta}$ is not real
	valued is obtained by letting $G$ be a cyclic group of order
	$n$ and taking $\eta$ to be the representation whose value on
	a fixed generator is $e^{\frac{2 \pi i}{n}}$.

	One can also define a chain with transition probability \[
	\frac{dim(\rho)}{dim(\eta) dim(\lambda)} <\chi^{\rho}
	\chi^{\lambda},\chi^{\eta}> \] which would not require $\eta$
	real valued in Lemma \ref{usereal} but this is less useful for the applications at hand, since then Proposition \ref{calcspec} below would fail as any reversible Markov chain has real eigenvalues.
	
Proposition \ref{relatechains} shows that $M_H$ is in fact a special case of $L_{\eta}$. 

\begin{prop} \label{relatechains} Let $M_H$ be the Markov chain on irreducible representations of $G$ corresponding to the choice of subgroup $H$. Let $L_{\eta}$ be the Markov chain on the irreducible representations of $G$ corresponding to the choice that $\eta$ is the representation of $G$ on cosets of $H$ (i.e. the induction of the trivial representation of $H$ to $G$). Then $M_H=L_{\eta}$. 
\end{prop}

\begin{proof} Throughout the proof we let $Res,Ind$ denote restriction and induction of characters.
\begin{eqnarray*}
L_{\eta}(\lambda,\rho) & = & \frac{|H|}{|G|} \frac{dim(\rho)}{dim(\lambda)}
<\chi^{\rho}, \chi^{\lambda} Ind_{H}^{G}[1]>_G\\ & = & \frac{|H|}{|G|}
\frac{dim(\rho)}{dim(\lambda)} <\chi^{\rho} \overline{\chi^{\lambda}},
Ind_{H}^{G}[1]>_G\\
 & = & \frac{|H|}{|G|} \frac{dim(\rho)}{dim(\lambda)} <Res_{H}(\chi^{\rho} \overline{\chi^{\lambda}}),1>_H \\
& = & \frac{|H|}{|G|} \frac{dim(\rho)}{dim(\lambda)} <Res_{H}(\chi^{\rho}), Res_H(\chi^{\lambda})>_H \\
& = & \frac{|H|}{|G|} \frac{dim(\rho)}{dim(\lambda)} \sum_{\tau \in Irr(H)} \kappa(\tau,\lambda) \kappa(\tau,\rho)\\ & = & M_H(\lambda,\rho). \end{eqnarray*} Note that the third equality is Frobenius reciprocity. \end{proof}

	Next we note that the chain $L_{\eta}$ can be explicitly
	diagonalized, a fact which has implications for the
	decomposition of tensor products. As this directly generalizes
	results from \cite{F2} (which explains their importance) and
	can be proved by a similar technique, we omit the proofs.

\begin{prop} \label{calcspec} Let $G$ be a finite group and $\eta$ any representation of $G$ whose character is real valued. Let $\pi$ denote the Plancherel measure of $G$. Then the eigenvalues and
 eigenfunctions of the Markov chain $L_{\eta}$ are indexed by conjugacy
 classes $C$ of $G$.
\begin{enumerate}
\item The eigenvalue parameterized by $C$ is
$\frac{\chi^{\eta}(C)}{dim(\eta)}$.

\item An orthonormal basis of eigenfunctions $\psi_C$ in $L^2(\pi)$ is defined by $\psi_C(\rho) = \frac{|C|^{\frac{1}{2}} \chi^{\rho}(C)}{dim(\rho)}$.
\end{enumerate}
\end{prop}

\begin{prop} Let $\eta$ be a representation of a finite group $G$ whose character $\chi^{\eta}$ is real valued. Suppose that $|G|>1$. Let $\beta \
= max_{g \neq 1} \frac{\chi^{\eta}(g)}{dim(\eta)}$ and let $\pi$
 denote the Plancherel measure of $G$.  Then for integer $r \geq 1$, \[ \sum_{\rho \in Irr(G)}
 |\frac{ dim(\rho)}{dim(\eta)^r} \ <\chi^{\rho}, (\chi^{\eta})^r > -
 \pi(\rho)| \leq |G|^{1/2} \beta^r.\] \end{prop}

\section{Properties of Jack Polynomials} \label{jackprop}

	The purpose of this section is to collect properties of Jack polynomials which will be crucial in the proof of Theorem \ref{maintheorem}. A thorough introduction to Jack polynomials is in Chapter 6 of \cite{Mac}. We conform to Macdonald's notation and let $J_{\lambda}^{(\alpha)}$ denote the Jack polynomial with parameter $\alpha$ associated to the partition $\lambda$. When $\alpha=1$, the Jack polynomials are Schur functions, and when $\alpha=2$ or $\alpha=\frac{1}{2}$, they are zonal polynomials corresponding to spherical functions of a Gelfand pair.

	As in the introduction, given a box $s$ in the diagram of
$\lambda$, let $a(s)$ and $l(s)$ denote the arm and leg of $s$
respectively. One defines quantities \[ c_{\lambda}(\alpha) = \prod_{s
\in \lambda} (\alpha a(s) + l(s) +1) \] \[ c_{\lambda}'(\alpha) =
\prod_{s \in \lambda} (\alpha a(s) + l(s) + \alpha) .\] Recall that
$m_i(\lambda)$ denotes the number of parts of $\lambda$ of size $i$
and that $l(\lambda)$ denotes the total number of parts of
$\lambda$. We let $z_{\lambda}=\prod_{i \geq 1} i^{m_i(\lambda)}
m_i(\lambda)!$, the size of the centralizer of a permutation of cycle
type $\lambda$ in the symmetric group.

	Let $\theta^{\lambda}_{\mu}(\alpha)$ denote the coefficient of
the power sum symmetric function $p_{\mu}$ in
$J_{\lambda}^{(\alpha)}$. Lemma \ref{jackorth} gives orthogonality
relations for these coefficients. We remark that when $\alpha=1$,
$\theta^{\lambda}_{\mu}(1)$ is equal to $\frac{n!}{z_{\mu}}
\frac{\chi^{\lambda}_{\mu}}{dim(\lambda)}$ where
$\chi^{\lambda}_{\mu}$ is the character value of the representation of
$S_n$ parameterized by $\lambda$ on elements of cycle type $\mu$. Thus
when $\alpha=1$, Lemma \ref{jackorth} specializes to the orthogonality
relations for characters of the symmetric group.

\begin{lemma} \label{jackorth} (\cite{Mac}, page 382)
\begin{enumerate}
\item \[ \sum_{|\mu|=n} z_{\mu} \alpha^{l(\mu)} \theta^{\rho}_{\mu}(\alpha) \theta^{\lambda}_{\mu}(\alpha) = \delta_{\rho,\lambda} c_{\rho}(\alpha) c'_{\rho}(\alpha).\]

\item \[ \sum_{|\rho|=n} \frac{\theta^{\rho}_{\mu}(\alpha)
\theta^{\rho}_{\nu}(\alpha)}{c_{\rho}(\alpha) c'_{\rho}(\alpha)} =
\delta_{\mu,\nu} \frac{1}{z_{\mu} \alpha^{l(\mu)}}.\] \end{enumerate}
\end{lemma}

	The following special values of
	$\theta^{\lambda}_{\mu}(\alpha)$ will be needed. 

\begin{lemma} \label{special} 
\begin{enumerate}
\item (\cite{Mac}, page 382) \[ \theta^{\lambda}_{(1^n)}(\alpha) = 1.\]
\item (\cite{Mac}, page 383) \[ \theta^{(n)}_{\mu}(\alpha) = \frac{n!}{z_{\mu}} \alpha^{n-l(\mu)}.\]
\item (\cite{Mac}, page 384) \[ \theta^{\lambda}_{(2,1^{n-2})}(\alpha) = n(\lambda') \alpha - n(\lambda) .\]
\item (\cite{Sta}, page 107) \[ \theta^{(n-1,1)}_{\mu}(\alpha) = \frac{\alpha^{n-l(\mu)} n!}{z_{\mu}} \frac{(\alpha(n-1)+1) m_1(\mu) -n}{\alpha n (n-1)}.\]
\end{enumerate}
\end{lemma} 

	Next we consider the ring of symmetric functions, with inner product defined by the orthogonality condition $<p_{\nu},p_{\mu}>_{\alpha}=\delta_{\nu,\mu} z_{\mu} \alpha^{l(\mu)}$. By Lemma \ref{jackorth}, this is equivalent to the condition that $<J_{\eta}^{(\alpha)},J_{\lambda}^{(\alpha)}>_{\alpha}= \delta_{\eta,\lambda} c_{\lambda}(\alpha)c_{\lambda}'(\alpha)$. For a symmetric function $f$, its adjoint $f^{\perp}$ is defined by the condition $<fg,h>_{\alpha}=<g, f^{\perp}h>_{\alpha}$ for all $g,h$ in the ring of symmetric functions. It is straightforward to check that $p_1^{\perp}= \alpha \frac{\partial}{\partial p_1}$ (for the case $\alpha=1$ see page 76 of \cite{Mac}).

	Let \[ \psi_{\lambda/\tau}'(\alpha) = \prod_{s \in
C_{\lambda/\tau}-R_{\lambda/\tau}} \frac{(\alpha a_{\lambda}(s) +
l_{\lambda}(s)+1)}{(\alpha a_{\lambda}(s) + l_{\lambda}(s)+\alpha)}
\frac{(\alpha a_{\tau}(s) + l_{\tau}(s)+\alpha)}{(\alpha a_{\tau}(s) +
l_{\tau}(s)+1)} \] where $C_{\lambda/\tau}$ is the union of columns of $\lambda$ that intersect $\lambda-\tau$ and $R_{\lambda/\tau}$ is the union of rows of $\lambda$ that intersect $\lambda-\tau$.

\begin{lemma} \label{useful} \[ p_1^{\perp} J_{\lambda}^{(\alpha)} = \sum_{|\tau|=n-1} \frac{c_{\lambda}'(\alpha) \psi'_{\lambda/\tau}(\alpha)}{c_{\tau}'(\alpha)} J_{\tau}^{(\alpha)}.\] 
\end{lemma}

\begin{proof} Take the inner product of both sides with $J_{\tau}^{(\alpha)}$. The left hand side becomes \[
<p_1^{\perp} J_{\lambda}^{(\alpha)},J_{\tau}^{(\alpha)}>_{\alpha} = <
J_{\lambda}^{(\alpha)}, p_1 J_{\tau}^{(\alpha)}>_{\alpha} .\] Using the
Pieri rule for Jack symmetric functions (\cite{Mac}, page 340), this
becomes
\[ \frac{c_{\tau}(\alpha)}{c_{\lambda}(\alpha)}
\psi'_{\lambda/\tau}(\alpha)
<J_{\lambda}^{(\alpha)},J_{\lambda}^{(\alpha)}>_{\alpha} =
c_{\tau}(\alpha) c_{\lambda}'(\alpha)
\psi'_{\lambda/\tau}(\alpha).\] By the orthogonality
relations for the $J's$, this is equal to the inner product of the
right hand side with $J_{\tau}^{(\alpha)}$. \end{proof}

\section{Construction of an Exchangeable Pair} \label{exchange}

	The purpose of this section is twofold. First, we use the
theory of harmonic functions on Bratelli diagrams to construct an
exchangeable pair $(W_{\alpha},W^*_{\alpha})$ with respect to
Jack$_{\alpha}$ measure on the set of partitions of size $n$ (and as
usual, we suppose without loss of generality that $\alpha \geq 1$). We
give a Markov chain $M_{\alpha}$ which is a deformation of the chain
$M_H$ from Section \ref{revisit} (when $\alpha=1$ it corresponds to
the case that $G=S_n$ and $H=S_{n-1}$). The second and more subtle
part of this section is to show that this construction is closely
related to a chain $L_{\alpha}$ which is a deformation of the chain
$L_{\eta}$ from Section \ref{revisit} (when $\alpha=1$ it corresponds
to the case that $G=S_n$ and $\eta$ is the irreducible representation
of the symmetric group of shape $(n-1,1)$). In fact much of this paper
can be pushed through for generalizations of $M_{\alpha}$ and
$L_{\alpha}$ corresponding to more vigorous walks on the set of
partitions, but for Stein's method it is preferable to use local
walks.

The use of both $M_{\alpha}$ and $L_{\alpha}$ will be crucial to this paper. An interesting result in this section will be that (except for holding probabilities), $L_{\alpha}$ is a rescaling of $M_{\alpha}$, so that one can work with whichever is more convenient. For instance it will be clear from the definition that the transition probabilities of $M_{\alpha}$ are always non-negative. But except for cases such as $\alpha=1,2$ where there is a group theoretic reason, it will not be clear that the transition probabilities of $L_{\alpha}$ are always non-negative. But to prove that $W_{\alpha}$ is an eigenvector of $M_{\alpha}$, it will be convenient to use connections with $L_{\alpha}$.

	In order to define $M_{\alpha}$, we first recall results on the theory of harmonic functions on Bratelli
diagrams. The basic language was reviewed in Section
\ref{revisit}. The level $\Gamma_n$ consists of all partitions of size $n$. The multiplicity function $\kappa_{\alpha}(\tau,\lambda)$ is
defined as $\psi'_{\lambda/\tau}(\alpha)$ where $\psi'_{\lambda/\tau}(\alpha)$ was defined in Section \ref{jackprop}. A
result of Stanley \cite{Sta} is that $dim_{\alpha}(\lambda) = \frac{n!
\alpha^n}{c_{\lambda}'(\alpha)}$. Then \cite{K3} shows that the
Jack$_\alpha$ measure \[ \pi_{\alpha}(\lambda) = \frac{\alpha^n
n!}{c_{\lambda}(\alpha) c'_{\lambda}(\alpha)} \] forms a coherent set of probability distributions for this Bratelli diagram.

	Motivated by the discussion in Section \ref{revisit}, for
$\lambda,\rho \in \Gamma_n$, we define (for $\alpha \geq 1$) the
transition probability $M_{\alpha}(\lambda,\rho)$ to be
\begin{eqnarray*} & & \frac{\pi_{\alpha}(\rho)}{dim_{\alpha}(\lambda)
dim_{\alpha}(\rho)} \sum_{|\tau|=n-1} \frac{dim_{\alpha}(\tau)^2
\kappa_{\alpha}(\tau,\rho)
\kappa_{\alpha}(\tau,\lambda)}{\pi_{\alpha}(\tau)}\\ & = &
\frac{c_{\lambda}'(\alpha) }{\alpha n c_{\rho}(\alpha) }
\sum_{|\tau|=n-1} \frac{\psi'_{\lambda/\tau}(\alpha) \psi'_{\rho/\tau}(\alpha)
c_{\tau}(\alpha)}{c_{\tau}'(\alpha)}. \end{eqnarray*} Note that this
corresponds to transitioning down a level and then up a level in the Bratelli diagram. The
expression for $M_{\alpha}(\lambda,\rho)$ is a mess, but three useful
observations can be made. First being a sum of non-negative terms, it
is nonnegative. Second, it is clear that the transition mechanism
$M_{\alpha}$ proceeds by local moves, in the sense that if
$M_{\alpha}(\lambda,\rho) \neq 0$, then $\lambda$ and $\rho$ have a
common descendant. Third, $M_{\alpha}$ is reversible with respect to Jack$_{\alpha}$ measure.

As an example, when $n=3$ the reader can verify that the $M_{\alpha}$ transition probabilities (rows add to 1) are
\[ \begin{array}{c c c c}
                   & {\bf (3)} & {\bf (2,1)} & {\bf (1^3)} \\
                {\bf (3)} & \frac{1}{2 \alpha+1} & \frac{2 \alpha}{2 \alpha+1} & 0\\
 			{\bf (2,1)} & \frac{\alpha+2}{3(\alpha+1)(2 \alpha+1)} & \frac{2(\alpha^2+ 7 \alpha+1)}{3(\alpha+2)(2 \alpha+1)} & \frac{\alpha(2 \alpha+1)}{3 (\alpha+1)(\alpha+2)} \\
                {\bf (1^3)} & 0 & \frac{2}{\alpha+2} & \frac{\alpha}{\alpha+2}
          \end{array} \] 

	Next we define (for $\alpha \geq 1$) a chain $L_{\alpha}$ to have transition
	``probability'' \[ L_{\alpha}(\lambda,\rho) =
	\frac{1}{c_{\rho}(\alpha) c_{\rho}'(\alpha) \alpha^n n!}
	\sum_{|\mu|=n} (z_{\mu})^2 \alpha^{2 l(\mu)} \theta^{\lambda}_{\mu}(\alpha)
	\theta^{\rho}_{\mu}(\alpha) \theta^{(n-1,1)}_{\mu}(\alpha)
	.\] As an example, when $n=3$ using the special values of the $\theta$'s given in Lemma \ref{special} (and also the value $\theta^{\lambda}_{(3)}(\alpha)$ which is determined from the other values by the orthogonality relations Lemma \ref{jackorth}), the reader can verify that the $L_{\alpha}$ transition probabilities (rows add to 1) are
\[ \begin{array}{c c c c}
                   & {\bf (3)} & {\bf (2,1)} & {\bf (1^3)} \\ {\bf
                (3)} & 0 & 1 & 0\\ {\bf (2,1)} & \frac{\alpha+2}{6
                \alpha(\alpha+1)} & \frac{2\alpha^2+ 11 \alpha-4}{6
                \alpha (\alpha+2)} & \frac{(2 \alpha+1)^2}{6
                (\alpha+1)(\alpha+2)} \\ {\bf (1^3)} & 0 &
                \frac{2\alpha +1}{\alpha(\alpha+2)} &
                \frac{\alpha^2-1}{\alpha(\alpha+2)} \end{array} \]
                Since the $\theta$'s can be negative it is not clear
                (see more discussion below) that these transition
                ``probabilities'' are non-negative. However
                $L_{\alpha}$ is clearly  ``reversible'' with respect to Jack$_{\alpha}$ measure. Proposition \ref{sum1} shows
                that the transition probabilities sum to one.

\begin{prop} \label{sum1} \[ \sum_{|\rho|=n} L_{\alpha}(\lambda,\rho) = 1.\] \end{prop}

\begin{proof} By definition $\sum_{|\rho|=n} L_{\alpha}(\lambda,\rho)$ is equal to \[ \sum_{|\rho|=n}
 \frac{1}{c_{\rho}(\alpha) c_{\rho}'(\alpha) \alpha^n n!}
 \sum_{|\mu|=n} (z_{\mu})^2 \alpha^{2 l(\mu)} \theta^{\lambda}_{\mu}(\alpha)
 \theta^{\rho}_{\mu}(\alpha) \theta^{(n-1,1)}_{\mu}(\alpha)
.\] Using the fact from part 1 of Lemma \ref{special} that $\theta^{\rho}_{(1^n)}=1$, this can be rewritten as 
\[ \sum_{|\mu|=n} (z_{\mu})^2 \alpha^{2 l(\mu)}
 \theta^{\lambda}_{\mu}(\alpha) \theta^{(n-1,1)}_{\mu}(\alpha)
 \sum_{|\rho|=n}
 \frac{\theta^{\rho}_{\mu}(\alpha) \theta^{\rho}_{(1^n)} }{c_{\rho}(\alpha) c_{\rho}'(\alpha) \alpha^n n!}
.\] The result now follows from part 2 of Lemma \ref{jackorth} and part 1 of Lemma \ref{special}. \end{proof}

	Theorem \ref{surprise} establishes a fundamental relationship between the chains $M_{\alpha}$ and $L_{\alpha}$.

\begin{theorem} \label{surprise} If $\lambda \neq \rho$, then \[ L_{\alpha}(\lambda,\rho) = \frac{\alpha(n-1)+1}{\alpha(n-1)} M_{\alpha}(\lambda,\rho).\]
\end{theorem}

\begin{proof} By part 4 of Lemma \ref{special}, $L_{\alpha}(\lambda,\rho)$ is equal to \[
 \frac{1}{\alpha n(n-1) c_{\rho}(\alpha) c_{\rho}'(\alpha)}
 \sum_{|\mu|=n} \theta^{\lambda}_{\mu}(\alpha)
 \theta^{\rho}_{\mu}(\alpha) \alpha^{l(\mu)} z_{\mu}
 \left( (\alpha(n-1)+1)m_1(\mu)-n \right).\] Since $\lambda \neq \rho$,
 part 1 of Lemma \ref{jackorth} shows that this is equal to \[
 \frac{(\alpha(n-1)+1)}{\alpha n(n-1) c_{\rho}(\alpha)
 c_{\rho}'(\alpha)} \sum_{|\mu|=n} \theta^{\lambda}_{\mu}(\alpha)
 \theta^{\rho}_{\mu}(\alpha) \alpha^{l(\mu)} z_{\mu} m_1(\mu).\]

	Bearing in mind the results from Section \ref{jackprop}, this
can be rewritten as \begin{eqnarray*} & &
\frac{(\alpha(n-1)+1)}{\alpha n(n-1) c_{\rho}(\alpha)
c_{\rho}'(\alpha)} \sum_{|\mu|=n} <p_1 \frac{\partial}{\partial p_1}
\sum_{|\mu|=n} \theta^{\lambda}_{\mu}(\alpha) p_{\mu}, \sum_{|\mu|=n}
\theta^{\rho}_{\mu}(\alpha) p_{\mu}>_{\alpha}\\ & = &
\frac{(\alpha(n-1)+1)}{\alpha^2 n(n-1) c_{\rho}(\alpha)
c_{\rho}'(\alpha)} <p_1^{\perp} J_{\lambda}^{(\alpha)}, p_1^{\perp}
J_{\rho}^{(\alpha)}>_{\alpha}\\ & = & \frac{(\alpha(n-1)+1)}{\alpha^2
n(n-1) c_{\rho}(\alpha) c_{\rho}'(\alpha)} \\
& & <\sum_{|\tau|=n-1}
\frac{\psi_{\lambda/\tau}'(\alpha)
c_{\lambda}'(\alpha)}{c_{\tau}'(\alpha)} J_{\tau}^{(\alpha)},
\sum_{|\tau|=n-1} \frac{\psi_{\rho/\tau}'(\alpha)
c_{\rho}'(\alpha)}{c_{\tau}'(\alpha)} J_{\tau}^{(\alpha)}>_{\alpha}\\
& = & \frac{(\alpha(n-1)+1)}{\alpha (n-1)} \sum_{|\tau|=n-1} \frac{c_{\lambda}'(\alpha)
c_{\tau}(\alpha) \psi_{\lambda/\tau}'(\alpha)
\psi_{\rho/\tau}'(\alpha)}{\alpha n c_{\tau}'(\alpha) c_{\rho}(\alpha)}\\ & = &
\frac{\alpha(n-1)+1}{\alpha(n-1)} M_{\alpha}(\lambda,\rho). \end{eqnarray*} \end{proof}

	Note that Theorem \ref{surprise} implies that
	$L_{\alpha}(\lambda,\rho) \geq 0$ for $\lambda \neq \rho$. We
	conjecture that $L_{\alpha}(\lambda,\lambda) \geq 0$ for all
	$\lambda$ and $\alpha \geq 1$. Using Theorem \ref{surprise}, this is equivalent to
	the assertion that $M_{\alpha}(\lambda,\lambda) \geq
	\frac{1}{\alpha(n-1)+1}$ for all $\lambda$. However as this
	paper only uses non-negativity of $M_{\alpha}$, this conjecture
	is somewhat of a distraction and we do not pursue it here. The
	proof should not be too difficult.

	In fact since $L_1(\lambda,\rho)$ is simply the chain $L_{\eta}$ of
	Section \ref{revisit} with $\eta$ the irreducible
	representation of shape $(n-1,1)$, nonnegativity of $L_1$ is
	clear. To conclude this section we give a similar group
	theoretic argument that $L_{2}(\lambda,\rho) \geq 0$ for all
	$\lambda,\rho$.

\begin{prop} $L_{2}(\lambda,\rho) \geq 0$ for all $\lambda,\rho$.
\end{prop}

\begin{proof} Let $H_{2n}$ be the hyperoctahedral group of order $2^nn!$. Using the
 notation of Section 7.2 of \cite{Mac}  for
 the Gelfand pair $(S_{2n},H_{2n})$, given $\lambda,\mu$ partitions
 of $n$, let $\omega^{\lambda}_{\mu}$ be the value of the spherical function $\omega^{\lambda}$ on a double coset of type $\mu$. It follows that \[ L_{2}(\lambda,\rho) = \frac{(2^n
 n!)^2}{c_{\rho}(2) c'_{\rho}(2)} \sum_{|\mu|=n} \frac{1}{2^{l(\mu)}
 z_{\mu}} \omega^{\lambda}_{\mu} \omega^{\rho}_{\mu}
 \omega^{(n-1,1)}_{\mu}.\] It is a general fact (page 396 of
 \cite{Mac}) that if $\omega_1,\cdots,\omega_t$ are spherical
 functions for a Gelfand pair $(G,K)$ and $a_{ij}^k$ are defined by \[
 \omega_i \omega_j = \sum_k a_{ij}^k \omega_k \] (where the
 multiplication $\omega_i \omega_j$ denotes the pointwise product)
 then $a_{ij}^k$ are real and $\geq 0$. The proposition now follows from the
 orthogonality relation \[ \sum_{\mu} \frac{1}{2^{l(\mu)} z_{\mu}}
 \omega^{\lambda}_{\mu} \omega^{\nu}_{\mu} = \delta_{\lambda,\nu}
 \frac{c_{\lambda}(2)c'_{\lambda}(2)}{(2^n n!)^2}\] on page 406 of
 \cite{Mac}.
\end{proof}

\section{Jack Polynomials and the Metropolis Algorithm} \label{metropjack}

	To begin we recall the Metropolis algorithm \cite{MRRTT} for
	sampling from a positive probability $\pi(x)$ on a finite set
	$X$. A marvelous survey of the Metropolis algorithm,
	containing references and many examples is \cite{DSa}. The
	Metropolis algorithm is especially useful when one can
	understand the ratios $r_{y,x}=\frac{\pi(y)}{\pi(x)}$, but can
	not easily compute $\pi(x)$ (for instance in Ising-type
	models). Let $S(x,y)$ (the base chain) be the transition
	matrix of a symmetric irreducible Markov chain on $X$. Define
	the Metropolis chain $T$ by letting $T(x,y)$, the probability
	of moving from $x$ to $y$ be defined by
\[ \left\{ \begin{array}{ll}
        S(x,y) r_{y,x} & \mbox{if \ $r_{y,x}<1$}\\
S(x,y) & \mbox{if \ $y \neq x$ \ and $r_{y,x} \geq 1$}\\
  S(x,x)+ \sum_{z \neq x \atop r_{z,x}<1} S(x,z) (1-r_{z,x})  & \mbox{if \ $y=x$ }
                                                \end{array}
                        \right.                  \] 

This chain has desirable properties. First, is easy to implement. From
$x$, pick $y$ with probability $S(x,y)$. If $y \neq x$ and $r_{y,x}
\geq 1$, the chain moves to $y$. If $y \neq x$ and $r_{y,x} <1$, flip
a coin with success probability $r_{y,x}$. If the coin toss succeeds,
the chain moves to $y$. Otherwise the chain stays at $x$. Second, the
chain $T(x,y)$ is irreducible and aperiodic with stationary
distribution $\pi$. Thus taking sufficiently many steps according to
the chain $T$ one obtains an arbitrarily good approximate sample of
$\pi$.

	A remarkable result of Hanlon \cite{Ha} relates the Metropolis
 algorithm to Jack symmetric functions. Fix $\alpha \geq 1$. Hanlon
 defines a Markov chain $T_{\alpha}$ on the symmetric group $S_n$ as
 follows. Let $\pi(x)$ be the probability measure on $S_n$ which
 chooses $x$ with probability proportional to $\alpha^{-c(x)}$ where
 $c(x)$ is the number of cycles of $x$ (ironically for sampling
 purposes one does not need to use the Metropolis algorithm as the
 constant of proportionality can be exactly computed in this
 case). Let $S(x,y)=\frac{1}{{n \choose 2}}$ if $x^{-1}y$ is a
 transposition, and 0 otherwise. Then Hanlon defines $T_{\alpha}(x,y)$
 to be the resulting Metropolis chain. To be explicit, if
 $\lambda_{x}$ is the partition whose rows are the cycle lengths of
 $x$, then the chance $T_{\alpha}(x,y)$ of moving from $x$ to $y$ is
 \[ \left\{ \begin{array}{ll} \frac{(\alpha-1) n(\lambda'_{x})}{\alpha
 {n \choose 2}} & \mbox{if \ $y=x$}\\ \frac{1}{{n \choose 2}} &
 \mbox{if \ $y=x(i,j)$ \ and $c(y)=c(x)-1$}\\ \frac{1}{\alpha {n
 \choose 2}} & \mbox{if \ $y=x(i,j)$ \ and $c(y)=c(x)+1$}\\ 0 &
 \mbox{otherwise} \end{array} \right.  \] Thus for $n=3$ the
 transition matrix is (rows sum to 1) \[ \begin{array}{c c c c c c c}
 & {\bf id} & {\bf (12)} & {\bf (13)} & {\bf (23)} & {\bf (123)} &
 {\bf (132)} \\ {\bf id} & 0 & \frac{1}{3} & \frac{1}{3} & \frac{1}{3}
 & 0 & 0 \\ {\bf (12)} & \frac{1}{3 \alpha} & \frac{\alpha-1}{3
 \alpha} & 0 & 0 & \frac{1}{3} & \frac{1}{3} \\ {\bf (13)} &
 \frac{1}{3 \alpha} & 0 & \frac{\alpha-1}{3 \alpha} & 0 & \frac{1}{3}
 & \frac{1}{3} \\ {\bf (23)} & \frac{1}{3 \alpha} & 0 & 0 &
 \frac{\alpha-1}{3 \alpha} & \frac{1}{3} & \frac{1}{3} \\ {\bf (123)}
 & 0 & \frac{1}{3 \alpha} & \frac{1}{3 \alpha} & \frac{1}{3 \alpha} &
 1-\frac{1}{\alpha} & 0 \\ {\bf (132)} & 0 & \frac{1}{3 \alpha} &
 \frac{1}{3 \alpha} & \frac{1}{3 \alpha} & 0 & 1-\frac{1}{\alpha}
 \end{array} .\]

	It is clear that the transition matrix for $T_{\alpha}$ commutes with the action of $S_n$ on itself by conjugation. Thus lumping the chain $T_{\alpha}$ to conjugacy classes gives a Markov chain on conjugacy classes of $S_n$. We denote this lumped Metropolis chain by $K_{\alpha}$. The transition probability $K_{\alpha}(\mu,\nu)$ is defined as $\sum T_{\alpha}(x,y)$ where $x$ is any permutation in the class $\mu$ and $y$ consists of all permutations in the class $\nu$. For instance when $n=3$ the transition matrix (rows sum to 1) is 
\[ \begin{array}{c c c c}
                   & {\bf (1^3)} & {\bf (2,1)} & {\bf (3)} \\
                {\bf (1^3)} & 0 & 1 & 0\\
 			{\bf (2,1)} & \frac{1}{3 \alpha } & \frac{\alpha-1}{3 \alpha} & \frac{2}{3} \\
                {\bf (3)} & 0 & \frac{1}{\alpha} & 1-\frac{1}{\alpha}
          \end{array} .\]

	Theorem \ref{metro} is due to Hanlon and is quite deep. In
	\cite{DHa} it is applied to analyze the convergence rate of the
	Metropolis chain $T_{\alpha}$. The case $\alpha=1$ of Theorem
	\ref{metro} is the usual Fourier analysis on the symmetric
	group (see \cite{DS} for details and an application to
	analyzing the convergence rate of random walk generated by
	random transpositions).

\begin{theorem} \label{metro} (\cite{Ha}) Suppose that $\alpha \geq 1$. Then the 
chance that the lumped Metropolis chain $K_{\alpha}$ on partitions moves from $(1^n)$ to the partition $\mu$ after $r$ steps is equal to \[ \alpha^n n!
\sum_{|\rho|=n} \frac{\theta^{\rho}_{\mu}(\alpha)}{c_{\rho}(\alpha) c_{\rho}'(\alpha)} \left( \frac{\alpha n(\rho') - n(\rho)}{\alpha {n \choose 2}} \right)^r.\] \end{theorem}

	The following consequence is worth recording.

\begin{cor} \label{mom} Suppose that $\alpha \geq 1$. Then the 
chance that the lumped Metropolis chain $K_{\alpha}$ on partitions of size $n$ moves from the partition $(1^n)$ to itself after $r$ steps is the $r$th moment of the statistic $\frac{W_{\alpha}}{\sqrt{\alpha {n \choose 2}}}$ under Jack$_{\alpha}$ measure. \end{cor} \begin{proof} By part 1 of Lemma \ref{special}, $\theta^{\rho}_{(1^n)}(\alpha)=1$. The result is now clear from Theorem \ref{metro}. \end{proof}

Corollary \ref{mom} allows one to compute the $r$th moment of
$W_{\alpha}$ in terms of return probabilities of the Metropolis chain
$K_{\alpha}$. This opens the door to the method of moments approach to
proving a central limit theorem for $W_{\alpha}$, as in \cite{Ho} for
the special case $\alpha=1$. However we prefer the Stein's method
approach, as it comes with an error term. But in passing we note a
consequence which indicates that the scaling of $W_{\alpha}$ has been
chosen correctly.

\begin{cor} \label{exvar} Suppose that $\alpha \geq 1$. Then $E(W_{\alpha})=0$, $E(W_{\alpha}^2)=1$, and $E(W_{\alpha}^3)=\frac{\alpha-1}{\sqrt{\alpha {n \choose 2}}}$. \end{cor} 

\begin{proof} The chance that $K_{\alpha}$ goes from $(1^n)$ to itself in one step is $0$. Hence $E(W_{\alpha})=0$. The chance that $K_{\alpha}$ goes from $(1^n)$ to itself in two steps is computed to be
 $\frac{1}{\alpha {n \choose 2}}$. Hence $E(W_{\alpha}^2) = 1$. The
 chance that $K_{\alpha}$ goes from $(1^n)$ to itself in three steps
 is equal to the chance of going from $(1^n)$ to $(2,1^{n-2})$ in two
 steps, and then back to $(1^n)$. This chance is
 $\frac{\alpha-1}{\alpha^2 {n \choose 2 }^2}$. Hence
 $E(W_{\alpha}^3)=\frac{\alpha-1}{\sqrt{\alpha {n \choose
 2}}}$. \end{proof}

\section{Central Limit Theorem for Jack Measure} \label{CLT}

	In this section we prove Theorem \ref{maintheorem}. Thus
	$\alpha \geq 1$ is fixed and we aim to show that
	$W_{\alpha}(\lambda)= \frac{\alpha n(\lambda') -
	n(\lambda)}{\sqrt{\alpha {n \choose 2}}}$ satisfies a central
	limit theorem when $\lambda$ is chosen from Jack$_{\alpha}$
	measure.

Let $(W_{\alpha},W^*_{\alpha})$ be the exchangeable pair constructed
 in Section \ref{exchange} using the Markov chain
 $M_{\alpha}$. Abusing notation due to possible negativity issues, it
 is also convenient to let $(W_{\alpha},W'_{\alpha})$ be the
 exchangeable pair constructed in Section \ref{exchange} using
 $L_{\alpha}$. To apply Stein's method it is necessary to work with
 the genuine exchangeable pair $(W_{\alpha},W^*_{\alpha})$, but
 Theorem \ref{surprise} will reduce computations involving it to the
 more tractable pair $(W_{\alpha},W'_{\alpha})$.

	Proposition \ref{1strelate} shows that the hypothesis needed
	to apply the Stein method bound (Theorem \ref{steinbound}) is
	satisfied. It also tells us that $W_{\alpha}$ is an
	eigenvector for the Markov chain $M_{\alpha}$, with eigenvalue
	$1-\frac{2}{n}$. It is perhaps unexpected that this eigenvalue
	is independent of $\alpha$.

\begin{prop} \label{1strelate} $E^{W_{\alpha}}(W_{\alpha}^*) 
= (1-\frac{2}{n}) W_{\alpha}$. \end{prop}

\begin{proof} Theorem \ref{surprise} implies that
\[ E^{\lambda}(W_{\alpha}^*-W_{\alpha}) = \frac{\alpha(n-1)}{\alpha(n-1)+1} E^{\lambda}(W_{\alpha}'-W_{\alpha}).\] Using the definition of the chain $L_{\alpha}$ and part 3 of Lemma \ref{special}, it follows that
\begin{eqnarray*} 
& & E^{\lambda}(W_{\alpha}')\\
 & = & \frac{1}{\sqrt{\alpha {n \choose 2}}} \sum_{|\rho|=n} L_{\alpha}(\lambda,\rho)
\theta^{\rho}_{(2,1^{n-2})}(\alpha) \\ 
& = & \frac{1}{\sqrt{\alpha {n \choose 2}}} \sum_{|\rho|=n}
\frac{\theta^{\rho}_{(2,1^{n-2})}(\alpha)}{c_{\rho}(\alpha) c_{\rho}'(\alpha) \alpha^n n!}
\sum_{|\mu|=n} (z_{\mu})^2
\alpha^{2 l(\mu)} \theta^{\lambda}_{\mu}(\alpha)
\theta^{\rho}_{\mu}(\alpha) \theta^{(n-1,1)}_{\mu}(\alpha) \\ & = &
\frac{1}{\sqrt{\alpha {n \choose 2}}} \sum_{|\mu|=n} (z_{\mu})^2 \alpha^{2 l(\mu)}
\theta^{\lambda}_{\mu}(\alpha) \theta^{(n-1,1)}_{\mu}(\alpha)
\sum_{|\rho|=n}
\frac{\theta^{\rho}_{\mu}(\alpha)
\theta^{\rho}_{(2,1^{n-2})}(\alpha)}{c_{\rho}(\alpha) c_{\rho}'(\alpha) \alpha^n n!}
. \end{eqnarray*}

Using part 2 of Lemma \ref{jackorth}, one sees that only the term
$\mu=(2,1^{n-2})$ makes a non-zero contribution. Thus
\begin{eqnarray*} E^{\lambda}(W_{\alpha}') & = & \frac{2}{\alpha
n(n-1)} \theta^{(n-1,1)}_{(2,1^{n-2})}(\alpha)
\frac{\theta^{\lambda}_{(2,1^{n-2})}(\alpha)}{\sqrt{\alpha {n \choose
2}}}\\ & = & \frac{2}{\alpha n(n-1)}
\theta^{(n-1,1)}_{(2,1^{n-2})}(\alpha) W_{\alpha}\\ & = & \left( 1-
\frac{2(\alpha n - \alpha+1)}{n(n-1) \alpha}
\right)W_{\alpha}. \end{eqnarray*} The last two equations used Lemma
\ref{special}. Consequently \[ E^{\lambda}(W_{\alpha}'-W_{\alpha}) =
-\left( \frac{2(\alpha n - \alpha+1)}{n(n-1) \alpha}
\right)W_{\alpha}. \] Thus
$E^{\lambda}(W_{\alpha}^*-W_{\alpha})=-\frac{2}{n} W_{\alpha}$, and
since this depends on $\lambda$ only through $W_{\alpha}$, the result
follows. \end{proof}

	More generally, the following proposition (proved using the same method as for Proposition \ref{1strelate}) holds.

\begin{prop} \label{notneed} Fix $\nu$ a partition of $n$. Then $\theta^{\lambda}_{\nu}(\alpha)$ is an eigenvector of $L_{\alpha}$ with eigenvalue $\frac{z_{\nu}}{\alpha^{n-l(\nu)} n!} \theta^{(n-1,1)}_{\nu}(\alpha)$ and an eigenvector of $M_{\alpha}$ with eigenvalue \[ 1+\frac{\alpha(n-1)}{\alpha(n-1)+1} \left( \frac{z_{\nu}}{\alpha^{n-l(\nu)} n!} \theta^{(n-1,1)}_{\nu}(\alpha)-1 \right).\] \end{prop} 

	As a consequence of Proposition \ref{1strelate}, we see that the mean $E(W_{\alpha})$ is equal to 0.

\begin{cor} \label{mean0} $E(W_{\alpha})=0$. \end{cor} \begin{proof} Since the pair $(W_{\alpha},W_{\alpha}^*)$ is 
exchangeable, $E(W_{\alpha}^*-W_{\alpha})=0$. Using Proposition \ref{1strelate}, we see that \[ 
E(W_{\alpha}^*-W_{\alpha})  =  E(E^{W_{\alpha}}(W_{\alpha}^*-W_{\alpha})) = -\frac{2}{n} E(W_{\alpha}) .\] Hence 
$E(W_{\alpha})=0$. \end{proof}

	Next we compute $E^{\lambda}(W_{\alpha}')^2$. Recall 
that this notation means the expected value of $(W_{\alpha}')^2$ given $\lambda$. This 
will be useful for analyzing the error term in Theorem \ref{steinbound}.  

\begin{prop} \label{forterm1} \begin{eqnarray*} E^{\lambda}((W_{\alpha}')^2)
& = & 1 + \theta^{\lambda}_{(2,1^{n-2})}(\alpha) \frac{4 (\alpha-1) (\alpha {n-1 \choose 2}-1)}{\alpha^2 n^2 (n-1)^2}\\
& & + \theta^{\lambda}_{(3,1^{n-3})}(\alpha) \frac{6 \left( \alpha (n-1)(n-3)-3 \right)}{\alpha^2 n^2 (n-1)^2}\\
& & + \theta^{\lambda}_{(2^2,1^{n-4})}(\alpha) \frac{4 \left(\alpha (n-1)(n-4)-4 \right)}{\alpha^2 n^2 (n-1)^2}. \end{eqnarray*}
\end{prop} 

\begin{proof} 
\begin{eqnarray*} E^{\lambda}((W'_{\alpha})^2) & = & \alpha {n \choose 2} \sum_{|\rho|=n} L_{\alpha}(\lambda,\rho) \left( \frac{\alpha n(\rho')-n(\rho)}{\alpha {n \choose 2}} \right)^2\\
& = & \alpha {n \choose 2} \sum_{|\rho|=n} \frac{1}{c_{\rho}(\alpha) c_{\rho}'(\alpha) \alpha^n n!}\\
& & \cdot \sum_{|\mu|=n} (z_{\mu})^2
\alpha^{2 l(\mu)} \theta^{\lambda}_{\mu}(\alpha)
\theta^{\rho}_{\mu}(\alpha) \theta^{(n-1,1)}_{\mu}(\alpha) \left( \frac{\alpha n(\rho')-n(\rho)}{\alpha {n \choose 2}} \right)^2\\
& = & \alpha {n \choose 2} \sum_{|\mu|=n} \theta^{\lambda}_{\mu}(\alpha) \theta^{(n-1,1)}_{\mu}(\alpha) \frac{(z_{\mu})^2 \alpha^{2 l(\mu)}}{\alpha^{n} n!}\\ & & \cdot \sum_{|\rho|=n} \frac{ \theta^{\rho}_{\mu}(\alpha)}{c_{\rho}(\alpha) c_{\rho}'(\alpha)} \left( \frac{\alpha n(\rho')-n(\rho)}{\alpha {n \choose 2}} \right)^2. \end{eqnarray*} 

	Next observe that using Theorem \ref{metro}, one can compute
	the sum \[ \alpha^n n! \sum_{|\rho|=n} \frac{
	\theta^{\rho}_{\mu}(\alpha)}{c_{\rho}(\alpha)
	c_{\rho}'(\alpha)} \left( \frac{\alpha
	n(\rho')-n(\rho)}{\alpha {n \choose 2}} \right)^2 \] for any
	partition $\mu$. Indeed, it is simply the probability that the
	lumped Metropolis chain $K_{\alpha}$ moves from $(1^n)$ to
	$\mu$ in two steps. From the explicit description of the
	transition rule of $K_{\alpha}$, it is straightforward to
	calculate that this probability is $\frac{1}{\alpha {n \choose
	2}}$ when $\mu=(1^n)$, is $\frac{\alpha-1}{\alpha {n \choose
	2}}$ when $\mu=(2,1^{n-2})$, is $\frac{4(n-2)}{n(n-1)}$ when
	$\mu=(3,1^{n-3})$, and is $\frac{(n-2)(n-3)}{n(n-1)}$ when
	$\mu=(2^2,1^{n-4})$. Together with part 4 of Lemma
	\ref{special}, this completes the proof of the
	proposition. \end{proof}

	One can use Proposition \ref{forterm1} to give a Stein's
	method proof of the fact that $Var(W_{\alpha})=1$, but in
	light of Corollary \ref{exvar} there is no need to do so.

	In order to prove Theorem \ref{maintheorem}, we have to
	analyze the error terms in Theorem \ref{steinbound}. To begin
	we study \[ E\left(-1+\frac{n}{4} E^{\lambda}(W_{\alpha}^*-W_{\alpha})^2
	\right)^2, \] obtaining an exact formula. From Jensen's inequality for conditional expectations, (see Lemma 5 of
	\cite{F4} for details) the fact that $W_{\alpha}$ is determined by
	$\lambda$ implies that \[ E[E^{W_{\alpha}}(W_{\alpha}^*-W_{\alpha})^2]^2 \leq
	E[E^{\lambda}(W_{\alpha}^*-W_{\alpha})^2]^2 .\] Hence Proposition \ref{term1}
	gives an upper bound on \[ E\left(-1+\frac{n}{4}
	E^{W_{\alpha}}(W_{\alpha}^*-W_{\alpha})^2 \right)^2. \]

\begin{prop} \label{term1} \[ E \left(-1+\frac{n}{4} E^{\lambda}(W_{\alpha}^*-W_{\alpha})^2 
\right)^2 = \frac{3 \alpha n + 2 \alpha^2-10 \alpha+2}{4 \alpha n(n-1)}.\] \end{prop}

\begin{proof} By Theorem \ref{surprise} and Proposition \ref{1strelate},
\begin{eqnarray*}
E^{\lambda}(W_{\alpha}^*-W_{\alpha})^2 & = & \frac{ \alpha(n-1)}{\alpha(n-1)+1} E^{\lambda}(W_{\alpha}'-W_{\alpha})^2\\
& = & \frac{ \alpha(n-1)}{\alpha(n-1)+1} \left(
W_{\alpha}^2 - 2W_{\alpha}E^{\lambda}(W_{\alpha}') + 
E^{\lambda}(W_{\alpha}')^2 \right)\\
& = & \frac{ \alpha(n-1)}{\alpha(n-1)+1}  \left( (\frac{4(\alpha n-\alpha+1)}{\alpha n(n-1)}-1)W_{\alpha}^2 + E^{\lambda}(W_{\alpha}')^2 \right). \end{eqnarray*} Combining this with Proposition \ref{forterm1}, it follows that $-1+\frac{n}{4} E^{\lambda}(W_{\alpha}^*-W_{\alpha})^2$ is equal to $A+B+C+D+E$ where  
\begin{enumerate}

\item $A= -1+\frac{n}{4} \frac{\alpha (n-1)}{\alpha(n-1)+1}$

\item $B= \frac{(\alpha-1)(\alpha {n-1 \choose 2}-1)}{\alpha n (n-1)(\alpha n - \alpha+1)} \theta^{\lambda}_{(2,1^{n-2})}(\alpha) $

\item $C= \frac{3(\alpha(n-1)(n-3)-3)}{2\alpha n (n-1) (\alpha n - \alpha+1)} \theta^{\lambda}_{(3,1^{n-3})}(\alpha)$

\item $D= \frac{\alpha(n-1)(n-4)-4}{\alpha n (n-1) (\alpha n - \alpha+1)} \theta^{\lambda}_{(2^2,1^{n-4})}(\alpha)$

\item \begin{eqnarray*} E & = & \frac{n}{4} \frac{\alpha (n-1)}{\alpha
n - \alpha+1} (\frac{4(\alpha n -\alpha+1)}{\alpha n(n-1)}-1) \alpha
{n \choose 2} \left(\frac{\alpha n(\lambda')-n(\lambda)}{\alpha {n
\choose 2}} \right)^2\\ & = & \frac{n}{4} \frac{\alpha (n-1)}{\alpha n
- \alpha+1} (\frac{4(\alpha n -\alpha+1)}{\alpha n(n-1)}-1)
\frac{1}{\alpha {n \choose 2}}
(\theta^{\lambda}_{(2,1^{n-2})}(\alpha))^2. \end{eqnarray*} \end{enumerate}

	We need to compute the Jack$_{\alpha}$ average of $(A+B+C+D+E)^2$. Since $A^2$ is 
constant, the average of $A^2$ is
$\left(-1+\frac{n}{4} \frac{\alpha (n-1)}{\alpha(n-1)+1} \right)^2$. The Jack$_{\alpha}$
averages of $B^2$,$C^2$,$D^2$ can all be computed using part 2 of Lemma \ref{jackorth}. 
To compute the Jack$_{\alpha}$ average of $E^2$ one uses Theorem \ref{metro} to reduce to 
computing the probability that after three steps taken by the chain
$K_{\alpha}$ started from the partition $(1^n)$, that one is at the partition
$(2,1^{n-2})$. From the description of the entries of the transition
matrix of $K_{\alpha}$, one computes this probability to be
$\frac{2(3\alpha n^2+\alpha n+ 2 \alpha^2-16 \alpha+2)}{\alpha^2 n^2 (n-1)^2}$. The Jack$_{\alpha}$ averages of $2AB$, $2AC$,
$2AD$, $2BC$, $2BD$, $2CD$ are all 0 by part 2 of Lemma \ref{jackorth}. The
Jack$_{\alpha}$ average of $2AE$ is computed using the second expression for $E$ and part 2 of Lemma
\ref{jackorth}. Finally, Theorem \ref{metro} reduces computation of
the Jack$_{\alpha}$ average of $2BE$ (respectively $2CE$ and $2DE$) to the
probability that after two steps taken by the chain $K_{\alpha}$ started at
$(1^n)$, that one is at the partition $(2,1^{n-2})$ (respectively
$(3,1^{n-3})$ and $(2^2,1^{n-4})$).  Thus all of the enumerations are
elementary and adding up the terms yields the proposition. \end{proof}

	The final ingredient needed to prove Theorem \ref{maintheorem} is an upper bound on $E|W^*-W|^3$. Typically this is the crudest term in applications of Stein's method. 

	Lemma \ref{crude} bounds the tail probabilities for $\lambda_1,\lambda_1'$ under Jack$_{\alpha}$ measure.

\begin{lemma} \label{crude} Suppose that $\alpha > 0$.
\begin{enumerate}
\item The Jack$_{\alpha}$ probability that $\lambda_1 \geq 2e \sqrt{\frac{n}{\alpha}}$ is at most $\frac{\alpha n^2}{4^{2e \sqrt{\frac{n}{\alpha}}}}$.
\item The Jack$_{\alpha}$ probability that $\lambda_1' \geq 2e \sqrt{\alpha n}$ is at most $\frac{n^2}{\alpha 4^{ 2e \sqrt{n \alpha}}}$.
\end{enumerate}
\end{lemma}

\begin{proof} Given a partition $\lambda$, let $\tau$ be the partition of $n-\lambda_1$
 given by removing the first row of $\lambda$. Then by the definition
 of Jack$_{\alpha}$ measure, it follows that the Jack$_{\alpha}$
 measure of $\lambda$ is at most $\frac{n!}{(n-\lambda_1)! \lambda_1!
 (\alpha(\lambda_1-1)+1) \cdots (\alpha+1)}$ multiplied by the
 Jack$_{\alpha}$ measure of $\tau$. It follows
 that the Jack$_{\alpha}$ probability that $\lambda_1=l$ is at most \[
 \frac{n!}{(n-l)! l!} \frac{1}{\alpha^{l-1} (l-1)!} \leq
 (\frac{n}{\alpha})^{l} \frac{\alpha l}{l!^2}.\] Using the inequality $y! \geq (y/e)^y$ and assuming that $l \geq 2e \sqrt{\frac{n}{\alpha}}$ this is at most
 \[ (\frac{ne^2}{\alpha l^2})^{l} \alpha l
 \leq \frac{\alpha n}{4^{2e \sqrt{\frac{n}{\alpha}}}}.\] The first
 assertion follows by summing over $l$ with $n \geq l \geq 2e\sqrt{\frac{n}{\alpha}}$.

	The second assertion follows from the first assertion by
symmetry. Indeed, since the Jack$_{\alpha}$ measure of $\lambda'$ is
the Jack$_{\frac{1}{\alpha}}$ measure of $\lambda$, the
Jack$_{\alpha}$ probability that $\lambda_1' \geq 2e \sqrt{\alpha n}$
is equal to the Jack$_{\frac{1}{\alpha}}$ probability that $\lambda_1
\geq 2e \sqrt{\alpha n}$. Now apply part 1 of the lemma with $\alpha$
replaced by $\frac{1}{\alpha}$. \end{proof}

\begin{prop} \label{term2} Suppose that $\alpha \geq 1$. Then there is a constant $C_{\alpha}$ depending on $\alpha$ such that \[ E|W^*-W|^3 \leq C_{\alpha} n^{-3/2}\] for all $n$.
\end{prop} 

\begin{proof} Recall that \[ W= \frac{1}{\sqrt{ \alpha {n \choose 2}}} (\alpha n(\lambda')-n(\lambda)).\] From the definition of $M_{\alpha}$, it is clear that $\lambda^*$ is obtained from $\lambda$ by removing a box from the diagram of $\lambda$ and reattaching it somewhere. It follows that \[|W^*-W| \leq \frac{1}{\sqrt{ \alpha {n \choose 2}}} (\alpha(\lambda_1+1) + \lambda_1'+1).\] Indeed, suppose that $\lambda^*$ is obtained from 
$\lambda$ by moving a box from row $a$ and column $b$ to a different row $c$ and
column $d$. Then \[ W^*-W = \frac{1}{\sqrt{ \alpha {n \choose 2}}}
\left( \alpha(\lambda_c-\lambda_a+1) + (\lambda_b' - \lambda_d'-1)
\right).\]

	Suppose that $\lambda_1 \leq 2e \sqrt{\frac{n}{\alpha}}$ and
that $\lambda_1' \leq 2e \sqrt{\alpha n}$. Then by the previous
paragraph \[ |W^*-W| \leq \frac{C_0}{\sqrt{n}} \] for a universal
constant $C_0$ (not even depending on $\alpha$). Note by the first paragraph, that even if $\lambda_1 >
2e \sqrt{\frac{n}{\alpha}}$ or $\lambda_1' > 2e \sqrt{\alpha n}$
occurs, then $|W^*-W| \leq C_1 \sqrt{\alpha}$ for a universal constant
$C_1$. The result now follows by Lemma \ref{crude}, which shows that these events occur with very low probability for $\alpha$ fixed. \end{proof}

	Summarizing, now we prove Theorem \ref{maintheorem} (the main result).

\begin{proof} We use Theorem \ref{steinbound} with the exchangeable pair $(W,W^*)$ constructed in Section \ref{exchange}. Proposition \ref{1strelate} shows this to be possible with 
$\tau=\frac{2}{n}$. The result now follows from Proposition
\ref{term1} (together with the paragraph before it) and Proposition \ref{term2}. \end{proof}

\section{Acknowledgements} The author was partially supported by National Security Agency grant MDA904-03-1-0049. We thank a referee for helpful comments.

\end{document}